\documentclass[12pt]{amsart}
\usepackage{somedefs,tracefnt,latexsym,,rawfonts,epsfig,amsfonts,amssymb,latexsy
m,enumerate,amscd}
\def\cal{\mathcal}
\def\Bbb{\mathbb}

\def\p{\partial}
\newtheorem{prop}{Proposition}[subsection]
\newtheorem{thm}{Theorem}[subsection]
\newtheorem*{lemma1}{Key Lemma}
\newtheorem{lemma}{Lemma}[subsection]
\newtheorem{cor}{Corollary}[subsection]
\newtheorem{defn}{Definition}[subsection]
\newtheorem{rem}{Remark}[subsection]

\begin{document}
\date{March 13, 2006.}
\title[Isomorphism Conjecture and $K$-theory]{The  
isomorphism conjecture for $3$-manifold groups and $K$-theory of
virtually poly-surface groups}
\author[S.K. Roushon]{Sayed K. Roushon}
\address{School of Mathematics\\
Tata Institute\\
Homi Bhabha Road\\
Mumbai 400005, India}
\email{roushon@math.tifr.res.in} 
\urladdr{http://www.math.tifr.res.in/\~\ roushon/paper.html}
\begin{abstract} 
This article has two purposes. In \cite{R3} we 
showed that the FIC (Fibered Isomorphism Conjecture for
pseudoisotopy functor) for a particular class of 
$3$-manifolds (we denoted this class by $\cal C$) is the key
to prove the FIC for $3$-manifold groups in general. And we  
proved the FIC for the fundamental groups of members of a subclass of
$\cal C$.
This result was obtained by showing that the double of any member of
this subclass is either Seifert fibered or supports a nonpositively curved
metric. In this article we 
prove that for any $M\in {\cal C}$ there is a closed $3$-manifold $P$ such
that either $P$ is
Seifert fibered or is a nonpositively curved $3$-manifold and $\pi_1(M)$
is a subgroup of $\pi_1(P)$. As a
consequence this proves that the FIC is true for any $B$-group
(see definition 4.2 in \cite{R3}). Therefore, the FIC is
true for any Haken $3$-manifold group and hence for any
$3$-manifold group (using the reduction theorem of \cite{R3})
provided we assume the Geometrization conjecture. The above result also
proves the FIC for a class of $4$-manifold groups (see \cite{R2}). 

The second aspect of this article is to relax a condition in
the definition of strongly poly-surface group (\cite{R1}) and define a new
class of groups (we call them {\it weak strongly poly-surface} groups).
Then using the above result we prove the FIC for any virtually weak
strongly poly-surface group. We also give a corrected proof
of the main lemma of \cite{R1}.\end{abstract} 



\maketitle

\section{Introduction}

Throughout this article by `FIC' we mean the Fibered Isomorphism
Conjecture
of Farrell and Jones (\cite{FJ}) corresponding to the stable topological
pseudoisotopy functor. The statement and the basic results related to
this conjecture and the basics needed on $3$-manifold topology have
already
appeared in several articles, for example see
\cite{FL}, \cite{R2} and \cite{R3}. We refer the reader to
these sources for the preliminaries.

In this article we always consider orientable $3$-manifolds. Note that the
FIC for the fundamental group of nonorientable $3$-manifolds follow easily
from the orientable ones (see remark 3.1 in \cite{R3}). 

\subsection{Geometry and the FIC for $3$-manifold groups}

In \cite{R3} we called a group $G$ a $B$-group if it
contains the fundamental group of a compact irreducible $3$-manifold with
all the boundary components of higher genus (that is, when genus $\geq
2$) and incompressible (nonempty) 
as a subgroup of finite index. We denoted by $\cal 
C$ the later class of $3$-manifolds. In the proof of the main
theorem of \cite{R3} we showed that the $B$-groups are the key class of
groups for
which the FIC should be proved to prove the FIC for $3$-manifold groups.
If $M\in {\cal C}$ and $M$ contains no essential annuli (see
definition 4.3 in \cite{R3}) then we showed in
theorem 4.3 of \cite{R3} that the FIC is true for a
$B$-group which
contains $\pi_1(M)$ as a subgroup of finite index. This result was
accomplished by showing that in such a situation $\pi_1(M)$ is
isomorphic to a subgroup of the fundamental group of a closed
$3$-manifold $N$ so that either $N$ is a Seifert fibered space or is
a nonpositively curved Riemannian manifold. In fact in this case $N$ was
the double of $M$. In this article we prove the following.

\begin{thm} \label{npc} Let $M\in {\cal C}$. Then 
there is a closed $3$-manifold $P$ so that $P$ is either a
Seifert fibered space or is a nonpositively curved $3$-manifold
and $\pi_1(M)$ is a subgroup of $\pi_1(P)$.\end{thm}

\begin{rem}{\rm Here we remark that if $M$ has at least one torus boundary
component then it follows from theorem 3.2 of \cite{L} that $\pi_1(M)$ is
a subgroup of a closed nonpositively curved $3$-manifold. This can be seen
by first taking the double $N$ of $M$ along higher genus ($\geq 2$)
boundary components and then taking the double $P$ of $N$ and applying
\cite{L}.}\end{rem} 

An important application to the above theorem is the following theorem.

Before we state the theorem recall from \cite{R3} that we called a group
$\Gamma$ to satisfy the FICwF if the FIC is true for $\Gamma\wr F$ for any
finite group $F$.

\begin{thm} \label{key} (Key Lemma) The FICwF is true for 
$B$-groups.\end{thm}

The above Theorem has several consequences as described in 
\cite{R3}. We recall some of these consequences here. 

\begin{cor} The hypothesis {\tt the FIC is true for
$B$-groups} can
be removed from the statements of any result stated in \cite{R3}.\end{cor}

Among several results the following results were proved in \cite{R3}
under the assumption that the FIC is true for $B$-groups. For the other
consequences we refer the reader to \cite{R3}. 

\begin{cor} \label{map} (Proposition \ref{fibering}) Let $M^3$ be a closed
$3$-manifold fibering over the circle. Then the FICwF is true for
$\pi_1(M)$.\end{cor}

Recall that the main lemma of \cite{R1} is a particular case of the above 
corollary. Hence we have completed the proof of the main lemma of
\cite{R1}.

Corollary \ref{map} also has application in proving the FIC for mapping
class groups of surfaces (see \cite{BPL}).

Following the convention of \cite{R3}, in the next corollary
by 
`$3$-manifold' we mean irreducible $3$-manifold with infinite fundamental
group. Note that in such a situation the $3$-manifold is aspherical.

\begin{cor} \label{compact} ([\cite{R3}, theorem 4.5]) Let $M$ be either 
a compact $3$-manifold with nonempty
boundary or a noncompact $3$-manifold. Then the FICwF is true for
$\pi_1(M)$.\end{cor}

\begin{cor} ([\cite{R3}, theorem 4.5]) The FICwF is true for the
fundamental group of any virtual Haken $3$-manifold.\end{cor}

\begin{cor} ([\cite{R3}, corollary 4.2]) Assume that the Geometrization 
conjecture is true. Then the FICwF is true for any $3$-manifold 
group.\end{cor}

\begin{rem} {\rm The Farrell-Jones Isomorphism Conjecture in the case
of algebraic $K$-theory functor for a class of $3$-manifold groups is
considered in a recent paper (\cite{BL}) by Bartels and  
L\"{u}ck. It is shown there that the $K$-theoretic assembly map in
the conjecture is rationally injective for these $3$-manifold
groups.}\end{rem}
 
Now recall that in \cite{R2} we proved the FIC for the fundamental group
of a class of $4$-manifolds under a certain assumption ({\it special}) and
we pointed out there that we do not need this hypothesis if the FIC were
true for $B$-groups. Hence an application of Theorem \ref{key} is that
we can remove this assumption from these results and hence we have the
following more general result. 

\begin{cor} The hypothesis {\tt special} can be removed from the
statements of any result in \cite{R2}. In particular, we have that the
FIC is true for the fundamental group of the following $4$-manifolds $M$.

\begin{itemize}
\item $M$ is a surface bundle over a surface. ([\cite{R2}, theorem 1.5]).
\item $M$ is a fiber bundle over the circle and $M$ has a complex 
structure. ([\cite{R2}, theorem 1.3]).
\item $M$ is a complex elliptic surface. ([\cite{R2}, theorem 1.2]).
\end{itemize}
\end{cor}

\subsection{The FIC for poly-surface groups}

Now we come to the application part of Theorem \ref{npc} for poly-surface
groups. In \cite{R1} we defined strongly poly-surface
groups. This definition was along the same lines as the definition 
of strongly poly-free
groups in \cite{AFR}. But the strongly poly-surface groups are more
general and difficult to tackle. We describe the main problem which
occur with strongly poly-surface groups. In strongly poly-free groups a
certain kind of 
$3$-manifolds appeared. These were mapping tori of compact surfaces with
nonempty boundary. Leeb's work (\cite{L}) was directly applicable here to
use nonpositively curved geometry on such $3$-manifolds. But in general a
mapping torus of a closed surface does not support such a metric. In fact
an arbitrary mapping torus of a closed surface need not even has a
nonpositively curved metric in the Alexandrov sense. On the
other hand in
the case of strongly poly-surface groups we allowed mapping tori
of any surface (closed, compact or non-compact with one end). These are a
large class of
$3$-manifolds and as we recalled above, the geometry of nonpositively
curvature is not applicable here directly. In \cite{R3} we developed some
new techniques to tackle this kind of situation. The main idea was to use
direct limit argument and reduce it to compact $3$-manifold with nonempty
boundary case. In \cite{R3} we dealt with compact $3$-manifolds which
has at least one incompressible torus boundary component. In Theorem
\ref{npc} we complete this program considering the other cases.  

In this article we remove the `non-compact surface
with one end' hypothesis from the definition of strongly poly-surface
groups allowing non-compact surfaces with any number of ends and
define the following class of groups. Then using Corollaries \ref{map}
and \ref{compact} we prove the FIC in this general situation. 

\begin{defn} {\rm A discrete group $\Gamma$ is called 
{\it weak strongly poly-surface} if there exists a finite filtration of
$\Gamma$ by subgroups: $1=\Gamma_0\subset \Gamma_1\subset \cdots \subset
\Gamma_n=\Gamma$ such that the following conditions are satisfied:

\begin{itemize}
\item $\Gamma_i$ is normal in $\Gamma$ for each $i$.

\item $\Gamma_{i+1}/\Gamma_i$ is isomorphic to the fundamental group
of a surface $F_i$ (say).

\item for each $\gamma\in \Gamma$ and $i$  there is a diffeomorphism
$f:F_i\to F_i$ such
that the induced automorphism $f_{\#}$ of $\pi_1(F_i)$ is equal to
$c_\gamma$ up to inner automorphism, where $c_\gamma$ is the 
automorphism
of $\Gamma_{i+1}/\Gamma_i\ \simeq \ \pi_1(F_i)$ induced by the conjugation
action on $\Gamma$ by $\gamma$.
\end{itemize}

In such a situation we say that the group $\Gamma$ has {\it rank} $\leq
n$}.\end{defn}

The only difference between strongly poly-surface group and weak strongly
poly-surface group is in the last condition, that is when $\pi_1(F_i)$ is
infinitely generated. In strongly poly-surface group the condition was
that in this case $F_i$ has `one end' and here $F_i$ can have any number
of ends. All other conditions are same. 

We prove the following.

\begin{thm} \label{poly} Let $\Delta$ be a virtually weak strongly 
poly-surface group. Then FIC is true for $\Delta$.\end{thm}

Here recall that given a class of group $\cal A$, a group $\Gamma$ is
called virtually-$\cal A$, if there is a normal subgroup $A\in {\cal A}$
of $\Gamma$ of finite index.

\medskip
\noindent
{\bf Acknowledgment.} The author would like to thank G.A. Swarup and D.
Long for some discussions suggesting the use of pseudo-Anosov
diffeomorphism for the gluing procedure in the proof of Theorem \ref{npc}.

\section{Proof of Theorem \ref{npc} and Theorem \ref{key}}
\subsection{Proof of Theorem \ref{npc}}

Let us first recall the Nielsen-Thurston classification of surface
diffeomorphisms. Let $S$
be a closed orientable surface of genus $\geq 2$ and let $f:S\to S$ be an
orientation preserving diffeomorphism. Then $f$ belongs to one of the
following classes of diffeomorphisms.

\begin{itemize}
\item there is a positive integer $n$ so that $f^n$ is isotopic to the
identity.
\item $f$ is isotopic to a reducible diffeomorphism. that is, there is
a finite class of mutually nonparallel, disjoint and
essential simple closed curves $\cal A$ on $S$ so that an isotopy of $f$
leaves $\cal A$ invariant.
\item $f$ is pseudo-Anosov. 
\end{itemize}

We refer the reader to the papers \cite{CB}, \cite{FLP} and \cite{P} for
details and definitions of the above classification of surface
diffeomorphisms. But we recall the following two well-known facts we need.

\noindent
{\bf Fact 1.} {\it If $f$ is pseudo-Anosov then $f^n$ is 
pseudo-Anosov for all nonzero integers $n$. Also any isotopy of $f$ is
pseudo-Anosov.}

\noindent
{\bf Fact 2.} {\it The class of pseudo-Anosov diffeomorphisms is disjoint
from the first two classes of diffeomorphisms in the above
classification.}

The following lemma is an useful application of the above
classification. We need this lemma crucially for the proof of Theorem
\ref{npc}.

\begin{lemma} \label{pseudo} Let $S$ be a closed orientable surface of
genus $\geq 2$ and
let $\cal A$ be a finite class of disjoint, mutually nonparallel and
essential
simple closed curves on $S$. Let $f:S\to S$ be a pseudo-Anosov
diffeomorphism. Then there
is an integer $n_0$ so that no member of $\cal A$ is parallel to any
member of $f^{n_0}({\cal A})$.\end{lemma}

\begin{proof} As $f:S\to S$ is pseudo-Anosov, by Fact 1,
$f^n$ is
also pseudo-Anosov for all nonzero integers $n$. 

Let ${\cal A}=\{c_1,c_2,\ldots ,c_k\}$. Fix $l\in \{1,2,\ldots , k\}$.
Then we claim that 
for every nonzero integer $n$, $f^n(c_l)$ is not parallel to $c_l$.
Since otherwise $f^n$ will be reducible, which is a contradiction. Hence 
clearly, for distinct nonzero integers $n$ and $m$, $f^n(c_l)$ is not
parallel to $f^m(c_l)$. Therefore, there is a positive integer $N_l$ so
that $f^s(c_l)$ is not parallel to any $c_i$ for $i=1,2,\ldots , k$ and
for $s > N_l$. 

Let $n_0=max\{N_1, N_2,\ldots , N_k\}$. This proves the Lemma. \end{proof}

\begin{proof}[Proof of Theorem \ref{npc}]

If $M$ contains no essential annulus (see definition 4.3 in \cite{R3}) 
then
let $P$ be the double of $M$. Assume $P$ is not Seifert fibered. Then by
lemma 6.2 of \cite{R3} $P$ is not a graph manifold and hence has a
hyperbolic piece in its' JSJT (Jaco-Shalen-Johannson and
Thurston) decomposition. Using \cite{L} we get that $P$ supports a
nonpositively curved Riemannian metric.

The proof is a bit involved and is the main feature of this article when
there are essential annuli embedded in $M$. For instance in such a 
situation just by taking double of $M$ will produce more incompressible
tori in the double of $M$. 

So assume that there are essential annuli embedded in $M$. The main idea
is to identify two copies of $M$ along their boundary components by some 
diffeomorphisms of surfaces so that the essential annuli do not match to
produce unwanted incompressible tori. The suitable candidates for this
purpose are the pseudo-Anosov diffeomorphisms. 

Let $\p_1,\cdots , \p_k$ be the components of $\p M$. Then $\p_i$ is a
closed orientable surface of genus $\geq 2$ for each $i$. Since $M$
contains properly embedded essential annuli, these annuli produces a
finite class of pairwise disjoint, nonparallel, essential simple closed
curves on the boundary. Using the remark (Remark 2.1.1) below choose the
maximal
(finite) class of such annuli and let $\cal A$ be the corresponding
finite class of simple closed curves on $\p =\cup_{i=1}^{i=k}\p_i$. Let
${\cal
A}_i={\cal A}\cap \p_i$. By Lemma \ref{pseudo} there are pseudo-Anosov
diffeomorphisms $f_i:\p_i\to \p_i$ so that no member of ${\cal A}_i$ is
parallel to any member of $f_i({\cal A}_i)$.

We need the following remark for the remaining proof.

\begin{rem} {\rm There are finitely many essential annuli embedded in $M$
which are unique up to isotopy.
This follows from the uniqueness of the JSJ decomposition of $M$ along
incompressible tori and essential annuli. (see the splitting
theorem of \cite{JS}, p. 157).}\end{rem}

\medskip

\centerline{\psfig{figure=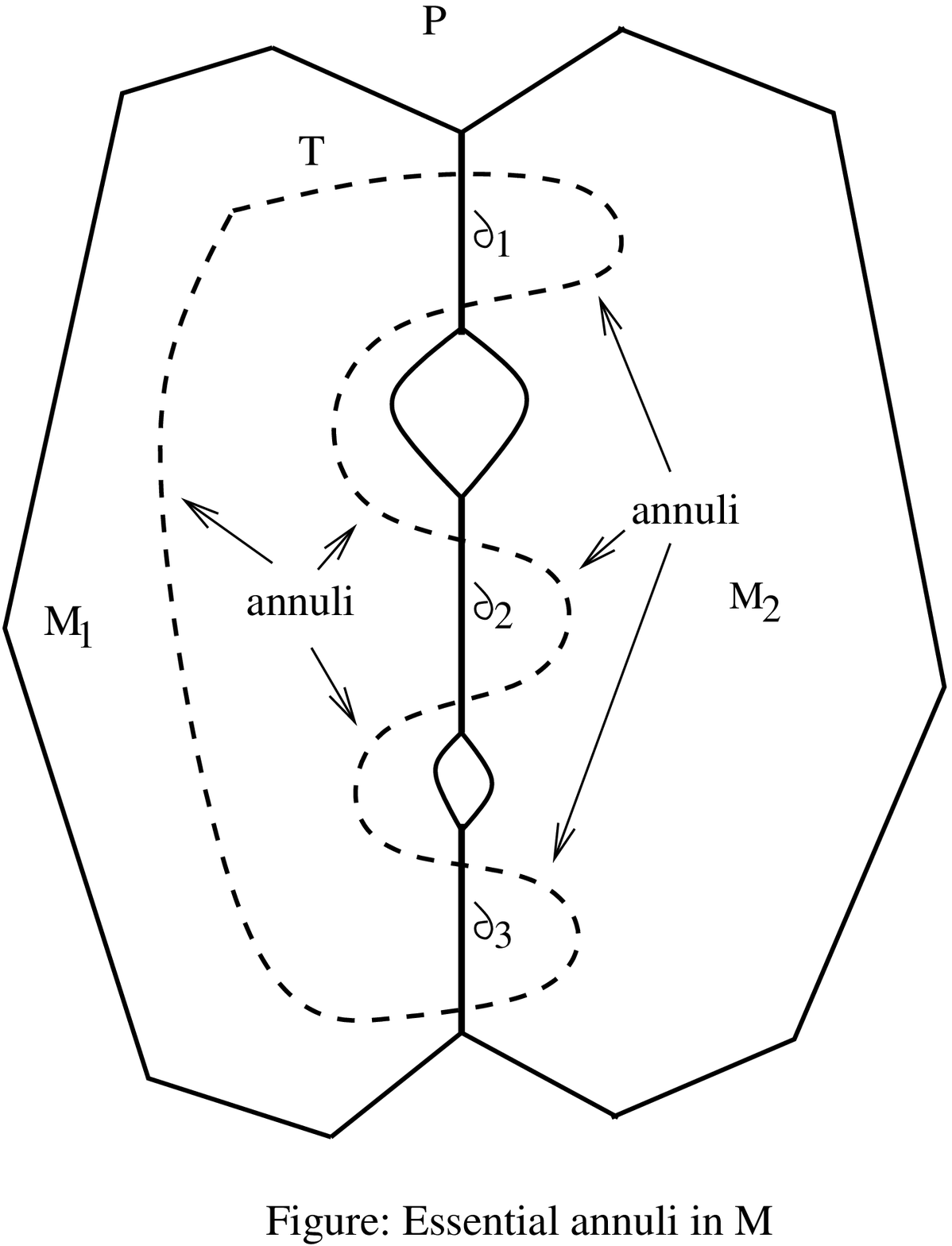,height=8cm,width=6cm}}

\vskip 0.5cm

Let $M_1$ and $M_2$ be two copies of $M$ and let us denote by ${\cal
A}^j$ the simple closed curves on $\p M_j$ for $j=1,2$ corresponding to
$\cal A$ on $\p M$. Let $P=M_1\cup M_2/\sim$, where $\sim$ stands for the
identification of $x$ with $f_i(x)$ for each $i$ and for
$x\in\p_i$. 

Note that $P$ is a Haken $3$-manifold and hence by the JSJT decomposition 
there are finitely many embedded incompressible tori (say, $\cal T$) in
$P$ so that the complementary pieces are either Seifert fibered or
hyperbolic. We check that either $P$ is Seifert fibered or it contains a
hyperbolic piece. This follows from the following two lemmas. 

Note that we can think of $\p$ as a submanifold of $P$. 

\begin{lemma} \label{iso} $\cal T$ can be isotoped off $\p$.\end{lemma}

\begin{proof} On the contrary, assume that there exists a $T\in {\cal T}$
such that
any isotopy of $T$ intersects $\p$. We can make this intersection in
general position
after isotoping $T$. Hence $T\cap M_j$ is a collection of annuli and
$2$-discs 
in $M_j$ for $j=1,2$. Now using the loop theorem (see \cite{He}, chapter
4, the loop theorem) and applying standard
tricks from $3$-manifold topology we can isotope $T$ further
to make sure
that the components of $T\cap M_j$ are essential annuli. Now by the
uniqueness of essential annuli
in $M$ it follows that $T\cap \p M_j$ is parallel to a subcollection of
${\cal A}^j$ for $j=1,2$. Since $T$ is obtained by gluing the
boundaries ($T\cap \p M_1$) of the annuli $T\cap M_1$ with the
boundaries ($T\cap \p M_2$) of the annuli $T\cap M_2$ under the maps
$\{f_i\}$, it follows that $f_i$ sends a member of ${\cal A}_i$
parallel to a member of $f_i({\cal A}_i)$. This is a contradiction to the
choice of $\{f_i\}$. Hence $\cal T$ can be isotoped off $\p$.\end{proof}

\begin{lemma} Either $P$ is Seifert fibered or there is a hyperbolic piece
in the JSJT decomposition of $P$.\end{lemma}

\begin{proof} Assume $P$ is not Seifert fibered. Let $S\in \{\p_1,\cdots ,
\p_k\}$. Using Lemma \ref{iso} isotope $\cal T$ so that it does not
intersect $\p=\cup_{i=1}^{i=k}\p_i$. Hence, there is a
component of $P-{\cal T}$, say $K$, containing $S$. We claim that $K$ is
not Seifert fibered. On the contrary, assume $K$ is Seifert fibered. Note
that $\pi_1(K)$ contains $\pi_1(S)$ as a
subgroup. Now since $P$ is not Seifert fibered, $K$ has nonempty boundary
and hence it has a
finite index subgroup of the form $F^r\times {\Bbb Z}$, where $F^r$ is
a (nonabelian) free group. Thus $\pi_1(S)\cap (F^r\times {\Bbb Z})$ is
of finite index in $\pi_1(S)$ and therefore is a closed surface group
of a surface of genus $\geq 2$. This is a contradiction because 
$F^r\times {\Bbb Z}$ cannot contain a nonabelian closed surface
group.\end{proof} 

Now using \cite{L} it follows in the non-Seifert fibered case that $P$
supports a nonpositively curved metric. This completes the proof of the
Theorem.\end{proof}

\subsection{A Key Lemma}
The key lemma to prove the FIC for discrete groups of this article is the
following and it is an application of 
Theorem \ref{npc}.

\begin{lemma1} Let $M\in \cal C$. Then the FICwF is true for $\pi_1(M)$.
Consequently, the FIC is true for  
$B$-groups.\end{lemma1}

\begin{rem} {\rm Here we recall that in the statement of the Key Lemma if
one of the boundary components of $M$ is an incompressible torus then the
same conclusion is
true. See corollary 4.1 of \cite{R3}}.\end{rem} 

The main application of the Key Lemma required for the proof of Theorem
\ref{poly} is the following proposition.

\begin{prop} \label{fibering} Let $M^3$ be a closed $3$-manifold fibering
over the circle. Then the FICwF is true for 
$\pi_1(M)$.\end{prop}

\begin{proof} The Key Lemma and theorem 4.5 of \cite{R3}  
together
completes the proof of the Proposition.\end{proof}

\begin{proof} [Proof of the Key Lemma] By Theorem \ref{npc} $\pi_1(M)$ is
isomorphic
to a subgroup of $\pi_1(P)$, where $P$ is either a Seifert fibered space 
or is a closed nonpositively curved
Riemannian manifold. Let $G$ be a finite group. Now it follows that
$\pi_1(M)\wr G$ is a
subgroup of $\pi_1(P)\wr G$. Then the FIC for $\pi_1(P)\wr G$ is a
consequence of theorem 4.1 of \cite{R3} when $P$ is
nonpositively
curved. The
Seifert fibered spaces case follows from theorem 4.6 of
\cite{R3}.
Using the hereditary property (lemma 2.1 of \cite{R3}) of
the FIC we complete the proof of the Key
Lemma.\end{proof} 

\section{Proof of Theorem \ref{poly}}
The proof goes along the same line as the proof of the main
theorem of \cite{R1}. 
The new ingredients in the following proof are Proposition \ref{fibering}
and Corollary \ref{compact}.

The proof is by induction on the rank of the finite index weak 
strongly poly-surface group.

\noindent
{\bf Induction hypothesis $I(n)$.} For any weak strongly poly-surface
group $\Gamma$ of rank $\leq n$ and for any finite group $G$, FIC is true
for the wreath product $\Gamma\wr G$.

If the rank of $\Gamma$ is $\leq0$ then $\Gamma\wr G=G$ finite and
hence $I(0)$ holds.

Now assume $I(n-1)$. We will show that $I(n)$ holds.  

Let $\Gamma$ be a weak strongly poly-surface group of rank $\leq n$ and is
a normal subgroup of $\Delta$ with $G$ as the finite quotient group. So we
have a filtration by subgroups $$1=\Gamma_0\subset \Gamma_1\subset \cdots
\subset \Gamma_n=\Gamma$$ with all the requirements as in the definition
of weak strongly poly-surface group. 

We have another exact sequence obtained after taking wreath 
product of the exact sequence $1\to \Gamma_1\to \Gamma\to
\Gamma/\Gamma_1\to 1$ with $G$. $$1\to
\Gamma_1^G\to \Gamma\wr G\to (\Gamma/\Gamma_1)\wr G\to 1$$ 

Let $p$ be the surjective homomorphism $\Gamma\wr G\to
(\Gamma/\Gamma_1)\wr G$. Note that $\Gamma/\Gamma_1$ is a weak strongly
poly-surface group of rank less or equal to $n-1$.

By induction hypothesis FIC is true for $(\Gamma/\Gamma_1)\wr G$. We
would like to apply lemma 2.2 of \cite{R3}. Let $Z$ be a
virtually cyclic
subgroup of $(\Gamma/\Gamma_1)\wr G$. Then there are two cases to
consider.

\noindent
{\bf $Z$ is finite.} In this case we have $p^{-1}(Z) < \Gamma_1^G\wr
Z < \Gamma_1\wr (G\times Z)$. Since $\Gamma_1$ is a surface group,
theorem 4.1 of \cite{R3} completes the proof in this case.

\noindent
{\bf $Z$ is infinite.} Let $Z_1=Z\cap (\Gamma/\Gamma_1)^G$. Then $Z_1$ is
an infinite cyclic normal subgroup of $Z$ of finite index. Let $Z_1$ be  
generated by $u$. We get $p^{-1}(Z) < p^{-1}(Z_1)\wr K$ where
$K$ is isomorphic to $Z/Z_1$.

Also, $$p^{-1}(Z_1)\wr K \ \simeq \ (\Gamma_1^G\rtimes \langle u\rangle
)\wr K < (\prod_{g\in G}(\Gamma_1\rtimes_{\alpha_g} \langle u\rangle ))\wr
K$$$$< \prod_{g\in G}((\Gamma_1\rtimes_{\alpha_g} \langle u\rangle ))\wr
K).$$

Now we describe the notations in the above display. Let $t\in \Gamma^G$
which goes to $u$. Then $\alpha_g(\gamma)=t_g\gamma t_g^{-1}$ for all
$\gamma\in \Gamma_1$ and $t_g$ is the value of $t$ at $g$. By definition  
each of these actions is induced by a diffeomorphism of the surface $F_0$
whose fundamental group is isomorphic to $\Gamma_1$.

Now there are two cases: (a) $\Gamma_1$ is finitely generated and (b)
$\Gamma_1$ is infinitely generated. 

The proof of Case $(a)$ follows from Proposition \ref{fibering} and
Corollary \ref{compact} and by noting that if the FIC is true for two
groups then it is true for their direct product also (see lemma 5.1 of
\cite{R3}). 

\noindent
{\bf (b).} $\Gamma_1$ is infinitely generated and hence free. We get that
each factor on the right hand side of the above display is a wreath
product of the fundamental group of a noncompact $3$-manifold with the
finite group $K$. Again using Corollary \ref{compact} we complete the
proof of the theorem.

\newpage 

\bibliographystyle{plain}
\ifx\undefined\bysame
\newcommand{\bysame}{\leavevmode\hbox to3em{\hrulefill}\,}
\fi

\end{document}